\documentclass[11pt]{article}
\usepackage[reqno]{amsmath}
\usepackage{amssymb}
\usepackage{amsfonts}

\newtheorem{theorem}{Theorem}

\newenvironment{proof}[1][Proof]{\noindent\textit{#1.} }{}
\input{tcilatex}

\begin{document}

\begin{center}
{\Huge On the High-Level Error Bound for Gaussian Interpolation}

\vspace*{0.2in}LIN-TIAN LUH

Department of Mathematics, Providence University, Shalu, Taichung, Taiwan

Abstract
\end{center}

\subparagraph{It's well-known that there is a very powerful error bound for
Gaussians put forward by Madych and Nelson in 1992. It's of the form$%
\left\vert f(x)-s(x)\right\vert \leq (Cd)^{\frac{c}{d}}\left\Vert
f\right\Vert _{h}$ where $C,c$ are constants, $h$ is the Gaussian function, $%
s$ is the interpolating function, and d is called fill distance which,
roughly speaking, measures the spacing of the points at which interpolation
occurs. This error bound gets small very fast as $d\rightarrow 0$. The
constants $C$ and $c$ are very sensitive. A slight change of them will
result in a huge change of the error bound. The number $c$ can be calculated
as shown in $\left[ 9\right] .$ However, $C$ cannot be calculated, or even
approximated. This is a famous question in the theory of radial basis
functions. The purpose of this paper is to answer this question.}


\begin{center}
{\LARGE 1.INTRODUCTION}
\end{center}

\subparagraph{Let $h$ be a continuous function on $R^{n\text{ }}$which is
conditionally positive definite of order $m.$ Given data $\left(
x_{j,}~f_{j}\right) ,j=1,...,N,$ where $X=\left\{ x_{1,}...,x_{N}\right\} $
is a subset of points in $R^{n}$ and the $f_{j}^{\prime }s$ are real or
complex numbers, the so-called $h$ spline interpolant of these data is the
function $s$ defined by}

\begin{flushright}
$s(x)=p(x)+\tsum\limits_{j=1}^{N}c_{j}h(x-x_{j}),$ \ \ \ \ \ \ \ \ \ \ \ \ \
\ \ \ \ \ \ \ \ \ \ \ \ \ $\left( 1\right) $
\end{flushright}

\subparagraph*{$\protect\bigskip $where $p(x)$ is a Polynomial in $\mathcal{P%
}_{m-1}$ and $c_{j}^{~~\prime }s$ are chosen so that}

\begin{flushright}
$\tsum\limits_{j=1}^{N}c_{j}q(x_{j})=0$ \ \ \ \ \ \ \ \ \ \ \ \ \ \ \ \ \ \
\ \ \ \ \ \ \ \ \ \ \ \ \ \ \ \ \ \ \ \ \ \ \ \ \ $\left( 2\right) $
\end{flushright}

\subparagraph{for all polynomials $q$ in $\mathcal{P}_{m-1}$ and$\ \ \ \ \ \
\ \ \ \ \ \ \ \ \ \ \ \ \ \ \ \ \ \ \ \ \ $}

\begin{flushright}
$p(x_{i})+\tsum\limits_{j=1}^{N}c_{j}h(x_{i}-x_{j})=f_{i},\quad i=1,...N.$ $%
\ \ \ \ \ \ \ \ \ \ \ \ \ \ \ \ \ \ \ \left( 3\right) $
\end{flushright}

\subparagraph{Here $\mathcal{P}_{m-1}$ denotes the class of those
polynomials of $R^{n}$ of degree $\leq m-1.$}

\subparagraph{It is well known that the system of equations $\left( 2\right) 
$ and $\left( 3\right) $ has a unique solution when $X$ is a determining set
for $\mathcal{P}_{m-1}$ and $h$ is strictly conditionally positive definite.
For more details see $\left[ 7\right] .$ Thus, in this case, the interpolant
s(x) is well defined.}

\subparagraph{We remind the reader that $X$ is said to be a determining set
for $\mathcal{P}_{m-1}$ if $p$ is in $\mathcal{P}_{m-1}$ and $p$ vanishes on 
$X$ implies that $p$ is identically zero.}

\subsection*{1.1A Bound for Multivariate Polynomials}

\subparagraph{A key ingredient in the development of our estimates is the
following lemma which gives a bound on the size of a polynimial on a cube in 
$R^{n}$ in terms of its values on a discrete subset which is scattered in a
sufficiently uniform manner. For its proof, please see $\left[ 9\right] $.}

\subparagraph{LEMMA1. For $n=1,2,...,$define $\protect\gamma _{n}$ by the
formulas $\protect\gamma _{1}=2$ and, if $n>1,$ $\protect\gamma _{n}=2n(1+%
\protect\gamma _{n-1}).$ Let $Q$ be a cube in $R^{n}$ that is subdivided
into $q^{n}$ identical subcubes. Let Y be a set of $q^{n}$ points obtained
by selecting a point from each of those subcubes. If $q\geq \protect\gamma %
_{n}(k+1),$ then for all $p$ in $\mathcal{P}_{k}$}

\begin{center}
$\underset{x\in Q}{\sup }\left\vert p(x)\right\vert \leq ~e^{2n\gamma
_{n}(k+1)}\underset{y\in Y}{\sup }\left\vert p(y)\right\vert .$
\end{center}

\subsection*{1.2 A Variational Framework for Interpolation}

\subparagraph{The precise statement of our estimate concerning $h$ splines
requires a certain amount of technical notation and terminology which is \
identical to that used in $\left[ 8\right] .$ For the convenience of the
reader we recall several basic notions.}

\subparagraph{The space of complex valued functions on $R^{n}$ that are
compactly supported and infinitely differentiable is denoted by $\mathcal{D}%
. $The Fourier transform of a function $\protect\phi $ in $\mathcal{D}$ is}

\begin{center}
$\overset{\wedge }{\phi }(\xi )=\int e^{-i<x,\xi >}\phi (x)dx.$
\end{center}

\subparagraph{In what follows $h$ will always denote a continuous
conditionally positive definite function of order $m.$ The Fourier transform
of such functions uniquely determines a positive Borel measure $\protect\mu $
on $R^{n}\thicksim \left\{ 0\right\} $ and constants $a_{\protect\gamma %
},\left\vert \protect\gamma \right\vert =2m$ as follows: For all $\protect%
\psi \in \mathcal{D}$}

\begin{center}
$\int h(x)\psi (x)dx=\dint \left\{ \overset{\wedge }{\psi }(\xi )-\overset{%
\wedge }{\chi }(\xi )\underset{\left\vert \gamma \right\vert <2m}{\sum }%
D^{\gamma }\overset{\wedge }{\psi }(0)\frac{\xi ^{\gamma }}{\gamma !}%
\right\} d\mu (\xi )$
\end{center}

\begin{flushright}
$\bigskip +\underset{\left\vert \gamma \right\vert \leq 2m}{\sum }D^{\gamma }%
\overset{\wedge }{\psi }(0)\frac{a_{\gamma }}{\gamma !},$ $\ \ \ \ \ \ \ \ \
\ \ \ \ \ \ \ \ \ \ \ \ \ \ \ \ \ \ \ \ \ \ \ \ \ \ \ \ \ \ \left( 4\right) $
\end{flushright}

\subparagraph{$\protect\bigskip $where for every choice of complex numbers $%
c_{\protect\alpha },\left\vert \protect\alpha \right\vert =m,$}

\begin{center}
$\underset{\left\vert \alpha \right\vert =m}{\sum }\underset{\left\vert
\beta \right\vert =m}{\sum }a_{\alpha +\beta }c_{\alpha }\overline{c_{\beta }%
}\geq 0.$
\end{center}

\subparagraph{Here $\protect\chi $ is a function in $\mathcal{D}$ such that $%
1-\protect\overset{\wedge }{\protect\chi }(\protect\xi )$ has a zreo of
order $2m+1$ at $\protect\xi =0;$ both of the integrals}

\begin{center}
$\int_{0<\left\vert \xi \right\vert <1}\left\vert \xi \right\vert ^{2m}d\mu
(\xi ),$ $\int_{\left\vert \xi \right\vert \geq 1}d\mu (\xi )$
\end{center}

\subparagraph{$\protect\bigskip $are finite. The choice of $\protect\chi $
affects the value of the coefficients $a_{\protect\gamma }$ for $\left\vert 
\protect\gamma \right\vert <$ $2m.$}

\subparagraph{Our variational framework for interpolation is supplied by a
space we denote by $\mathcal{C}_{h,m}.$ If}

\begin{center}
$\mathcal{D}_{m}=\left\{ \phi \in \mathcal{D}:\int x^{\alpha }\phi (x)dx=0%
\text{ }for~all~\left\vert \alpha \right\vert <m\right\} $
\end{center}

\subparagraph{,then $\mathcal{C}_{h,m}$ is the class of those continuous
functions $f$ which satisfy}

\begin{flushright}
$\left\vert \int f(x)\phi (x)dx\right\vert \leq c(f)\left\{ \int h(x-y)\phi
(x)\overline{\phi (y)}dxdy\right\} ^{\frac{1}{2}}\quad \ \ \ \ \ \ \ \ \ \ \
\ \ \ \ \left( 5\right) $
\end{flushright}

\subparagraph{for some constant $c(f)$ and all $\protect\phi $ in $\mathcal{D%
}_{m}.$ If $f\in \mathcal{C}_{h,m}$ let $\left\Vert f\right\Vert _{h}$
denote the smallest constant $c(f)$ for which $\left( 4\right) $ is true.
Recall that $\left\Vert f\right\Vert _{h}$ is a semi-norm and $\mathcal{C}%
_{h,m}$ is a semi-Hilbert space; in the case $m=0$ it is a norm and a
Hilbert space respectively.}

\subparagraph{Given a function $f$ in $\mathcal{C}_{h,m}$ and a subset X of $%
R^{n}$, there is an element $s$ of minimal $\mathcal{C}_{h,m}$ norm which is
equal to $f$ on $X.$}

\subparagraph{The function space $\mathcal{C}_{h,m}$ is not very easy to
understand. However Luh made a lucid characterization for this sapce in $%
\left[ 3\right] $ and $\left[ 4\right] $. For an easier understanding of
this space we suggest that reader read $\left[ 3\right] $ and $\left[ 4%
\right] $ first.}

\begin{center}
\bigskip {\Large 2.MAIN RESULTS}
\end{center}

\subparagraph{Before showing our main results, we need some lemmas. The
following lemma is cited directly from $\left[ 9\right] .$}

\subparagraph{LEMMA2. Let $Q,Y,~$and $\protect\gamma _{n}$ be as in LEMMA1.
Then, given a point $x$ in $Q,$ there is a measure $\protect\sigma $
supported on $Y$ such that}

\begin{center}
$\int p(y)d\sigma (y)=p(x)$
\end{center}

\subparagraph{for all $p$ in $\mathcal{P}_{k},$ and}

\begin{center}
$\int d\left\vert \sigma \right\vert (y)\leq e^{2n\gamma _{n}(k+1)}.$
\end{center}

\subparagraph{Now we need a famous formula.}

\subparagraph{\textbf{Stirling's Formula: }$n!\backsim \protect\sqrt{2%
\protect\pi n}(\frac{n}{e})^{n}.$\protect\smallskip}

\subparagraph{\protect\smallskip Remark}

\textit{The approximation is very reliable even for small }$n.$\textit{\ For
example, when }$n=10,$\textit{\ the relative error is only 0.83\%. The
larger }$n$\textit{\ is, the better the approximation is. For further
details see }$\left[ 1\right] $\textit{\ and }$\left[ 2\right] .$

\subparagraph{LEMMA3. Let $\protect\rho _{1}=\frac{1}{e}$ and $\protect\rho %
_{2}=\frac{3^{\frac{1}{6}}}{e}\backsim \frac{1.2}{e}$ . Then}

\begin{center}
$\sqrt{2\pi }\rho _{1}^{k}k^{k}\leq k!\leq \sqrt{2\pi }\rho _{2}^{k}k^{k}$
\end{center}

\textbf{for all positive integer} $k.$

\begin{proof}
Note that $\frac{1}{e},\frac{\sqrt{2}}{e^{2}},\frac{\sqrt{3}}{e^{3}},\frac{%
\sqrt{4}}{e^{4}},\frac{\sqrt{5}}{e^{5}},...$can be expressed by$\frac{1}{e},(%
\frac{2^{\frac{1}{4}}}{e})^{2},(\frac{3^{\frac{1}{6}}}{e})^{3},(\frac{4^{%
\frac{1}{8}}}{e})^{4},(\frac{5^{\frac{1}{10}}}{e})^{5},....~\qquad \qquad $

Now, $\sup \left\{ \frac{1}{e},\frac{2^{\frac{1}{4}}}{e},\frac{3^{\frac{1}{6}%
}}{e},\frac{4^{\frac{1}{8}}}{e},\frac{5^{\frac{1}{10}}}{e},...\right\} =%
\frac{3^{\frac{1}{6}}}{e}$implies that $\frac{\sqrt{k}}{e^{k}}\leq \rho
_{2}^{k}$ for all $k.$ Thus$k!\backsim \sqrt{2\pi }\frac{\sqrt{k}}{e^{k}}%
\cdot k^{k}\leq \sqrt{2\pi }\rho _{2}^{k}\cdot k^{k}.$The remaining part $%
\sqrt{2\pi }\rho _{1}^{k}k^{k}\leq k!$ follows by observing that $\sqrt{2\pi 
}(\frac{1}{e})^{k}k^{k}\leq \sqrt{2\pi }(\frac{1}{e})^{k}\cdot \sqrt{k}\cdot
k^{k}\backsim k!$

\begin{flushright}
$\square $
\end{flushright}
\end{proof}

\subparagraph{\textbf{LEMMA4. Let }$\protect\rho =\frac{\protect\sqrt{3}}{e}%
. $\textbf{Then }$k!\leq \protect\sqrt{2\protect\pi }\protect\rho ^{k}k^{k-1}
$\textbf{\ for all }$k\geq 1$\textbf{. \ \ \ \ \ \ \ \ \ \ \ }}

\begin{proof}
$k!\backsim \sqrt{2\pi }(\frac{1}{e})^{k}\cdot \sqrt{k}\cdot k^{k}=\sqrt{%
2\pi }\cdot \frac{k^{\frac{3}{2}}}{e^{k}}\cdot k^{k-1}.$ Note$\ $that$%
\left\{ \frac{k^{\frac{3}{2}}}{e^{k}}:k=1,2,3...\right\} $ can be expressed
by $\left\{ \frac{1}{e},(\frac{2^{\frac{3}{4}}}{e})^{2},(\frac{3^{\frac{1}{2}%
}}{e})^{3},(\frac{4^{\frac{3}{8}}}{e})^{4}....\right\} $ Our lemma follows
by noting that $\sup \left\{ \frac{1}{e},\frac{2^{\frac{3}{4}}}{e},\frac{3^{%
\frac{1}{2}}}{e},\frac{4^{\frac{3}{8}}}{e}...\right\} =\frac{\sqrt{3}}{e}.$
\end{proof}

\begin{flushright}
$\square $\bigskip
\end{flushright}

\subparagraph{LEMMA5. Let $h(x)=e^{-\protect\beta \left\vert x\right\vert
^{2}},\protect\beta >0,$ be the Gaussian function in $R^{n},$ and $\protect%
\mu $ be the measure defined in $\left( 4\right) .$ For any positive even
integer $k,$}

\begin{center}
$\tint\limits_{R^{n}}\left\vert \xi \right\vert ^{k}d\mu (\xi )\leq \pi ^{%
\frac{n+1}{2}}\cdot n\cdot \alpha _{n}\cdot 2^{\frac{k+n+2}{2}}\cdot \rho ^{%
\frac{k+n-1}{2}}\cdot \beta ^{\frac{k}{2}}\cdot (k+n-1)^{\frac{k+n-3}{2}%
}\cdot (2+\frac{1}{e})$
\end{center}

\subparagraph{\protect\bigskip for add n, and}

\begin{center}
$\tint\limits_{R^{n}}\left\vert \xi \right\vert ^{k}d\mu (\xi )\leq \pi ^{%
\frac{n+1}{2}}\cdot n\cdot \alpha _{n}\cdot 2^{\frac{k+n+3}{2}}\cdot \rho ^{%
\frac{k+n-2}{2}}\cdot \beta ^{\frac{k}{2}}\cdot (k+n-2)^{\frac{k+n-4}{2}}$
\end{center}

\subparagraph{for even n,where $\protect\rho =\frac{\protect\sqrt{3}}{e}$
and $\protect\alpha _{n}$ is the volume of the unit ball in $R^{n}.$}

\begin{proof}
$\tint\limits_{R^{n}}\left\vert \xi \right\vert ^{k}d\mu (\xi )$

$=(\frac{\pi }{\beta })^{\frac{n}{2}}\tint\limits_{R^{n}}\frac{\left\vert
\xi \right\vert ^{k}}{e^{\frac{\left\vert \xi \right\vert ^{2}}{4\beta }}}%
d\xi $

$=(\frac{\pi }{\beta })^{\frac{n}{2}}\cdot n\cdot \alpha _{n}\cdot
\int_{0}^{\infty }\frac{r^{k}\cdot r^{n-1}}{e^{\frac{\gamma ^{2}}{4\beta }}}%
d\gamma $

$=(\frac{\pi }{\beta })^{\frac{n}{2}}\cdot n\cdot \alpha _{n}\cdot (2\sqrt{%
\beta })^{k+n}\cdot \int_{0}^{\infty }\frac{r^{k+n-1}}{e^{r^{2}}}dr$

$=\left\{ 
\begin{array}{c}
\pi ^{\frac{n}{2}}\cdot n\cdot \alpha _{n}\cdot 2^{n-1+k}\cdot \beta ^{\frac{%
k}{2}}\cdot (\frac{k+n-2}{2})(\frac{k+n-2}{2}-1)\cdot \cdot \cdot (\frac{1}{2%
})\int_{0}^{\infty }\frac{1}{\sqrt{u}e^{u}}du\text{ }\ \text{if }n\text{ is
odd,\ \ \ \ \ \ \ \ \ \ \ \ \ \ \ \ \ \ \ \ \ \ \ \ \ \ \ \ \ \ \ \ \ \ \ \
\ \ \ \ \ \ \ \ \ \ } \\ 
\pi ^{\frac{n}{2}}\cdot n\cdot \alpha _{n}\cdot 2^{n-1+k}\cdot \beta ^{\frac{%
k}{2}}(\frac{k+n-2}{2})!\text{ if }n\text{ is even.\ \ \ \ \ \ \ \ \ \ \ \ \
\ \ \ \ \ \ \ \ \ \ \ \ \ \ \ \ \ \ \ \ \ \ \ \ \ \ \ \ \ \ \ \ \ \ \ \ \ \
\ \ \ \ \ \ \ \ \ \ \ \ \ \ \ \ \ \ \ \ \ \ \ \ \ \ \ \ \ \ \ \ \ \ \ \ \ \
\ \ \ }%
\end{array}%
\right. $
\end{proof}

\ \ \ \ \ \ \ \ \ \ \ \ \ \ \ $\ \ \ \ \ \ \ \ \ \ \ \ \ \ \ \ \ \ \ \ \ \ \
\ \ \ \ \ \ $

\subparagraph{Let $b=\left\lceil \frac{k\text{ }+n-2}{2}\right\rceil .$ Then 
$b!\leq \protect\sqrt{2\protect\pi }\protect\rho ^{b}\cdot b^{b-1}.$}

\subparagraph{For odd $n,$ $b=\frac{k+n-1}{2}.$ Thus}

\begin{center}
$\ \ \ \ \ \ \ \ \ \ \ \ \ \ \ b!\leq \sqrt{2\pi }\rho ^{\frac{k+n-1}{2}%
}\cdot (\frac{k+n-1}{2})^{\frac{k+n-3}{2}}$

$\ \ \ \ \ \ \ \ \ \ \ \ \ \ \ \ \ \ \ \ \ \ \ \ \ \ \ \ \ \ \ \ \ \ \ \ \ \
\ \ =\sqrt{2\pi }\cdot 2^{\frac{-k-n+3}{2}}\cdot \rho ^{\frac{k+n-1}{2}%
}\cdot (k+n-1)^{\frac{k+n-3}{2}},$
\end{center}

\subparagraph{and}

\begin{center}
$\tint\limits_{R^{n}}\left\vert \xi \right\vert ^{k}d\mu (\xi )\leq \pi ^{%
\frac{n+1}{2}}\cdot n\cdot \alpha _{n}\cdot (2+\frac{1}{e})\cdot 2^{\frac{%
k+n+2}{2}}\cdot \beta ^{\frac{k}{2}}\cdot \rho ^{\frac{k+n-1}{2}}\cdot
(k+n-1)^{\frac{k+n-3}{2}},$
\end{center}

\subparagraph{by noting that $\protect\int_{0}^{\infty }\frac{1}{\protect%
\sqrt{u}e^{u}}du\leq 2+\frac{1}{e}.$}

\subparagraph{For even $n,$ $b=\frac{k\text{ }+n-2}{2}.$ Thus}

\begin{center}
$b!\leq \sqrt{2\pi }\rho ^{\frac{k+n-2}{2}}\cdot (\frac{k+n-2}{2})^{\frac{%
k+n-4}{2}}$

$\ \ \ \ \ \ \ \ \ \ \ \ \ \ \ \ \ \ \ \ \ \ \ \ =\sqrt{2\pi }\cdot 2^{\frac{%
-k-n+4}{2}}\cdot \rho ^{\frac{k+n-2}{2}}\cdot (k+n-2)^{\frac{k+n-4}{2}},$
\end{center}

\subparagraph{and}

\begin{center}
$\tint\limits_{R^{n}}\left\vert \xi \right\vert ^{k}d\mu (\xi )\leq \pi ^{%
\frac{n+1}{2}}\cdot n\cdot \alpha _{n}\cdot 2^{\frac{k+n+3}{2}}\cdot \beta ^{%
\frac{k}{2}}\cdot \rho ^{\frac{k+n-2}{2}}\cdot (k+n-2)^{\frac{k+n-4}{2}}.$
\end{center}

\begin{flushright}
$\square $
\end{flushright}

\begin{theorem}
Let $h(x)=e^{-\beta \left\vert x\right\vert ^{2}},$ $\beta >0,$ be the
Gaussian function in $R^{n},$ and $\mu $ be the measure defined in $\left(
4\right) .$ Then, given a positive number $b_{0},$ there are positive
constants $\delta _{0},c,$and $C$ for which the following is true: If $f\in 
\mathcal{C}_{h,m}$ and $s$ is the $h$ spline that interpolates $f$ on a
subset $X$ of $R^{n},$ then
\end{theorem}

\begin{center}
$\left\vert f(x)-s(x)\right\vert \leq \triangle ^{^{\prime \prime
}}(C_{\delta })^{\frac{c}{\delta }}\cdot \left\Vert f\right\Vert _{h}$
\end{center}

\subparagraph{\protect\bigskip ,where $\triangle ^{^{\prime \prime }}=%
\protect\pi ^{\frac{n-1}{4}}\cdot (n\cdot \protect\alpha _{n})^{\frac{1}{2}%
}\cdot 2^{\frac{n+1}{4}}(\frac{\protect\sqrt{3}}{e})^{\frac{n-2}{4}}$ for
even $n$ and $\triangle ^{^{\prime \prime }}=\protect\pi ^{\frac{n-1}{4}%
}\cdot (n\cdot \protect\alpha _{n})^{\frac{1}{2}}\cdot (2+\frac{1}{e})^{%
\frac{1}{2}}\cdot 2^{\frac{n}{4}}\cdot (\frac{\protect\sqrt{3}}{e})^{\frac{%
n-1}{4}}$ for add $n,$ holds for all $x$ in a cube $E$ provided that $\left(
i\right) E$ has side $b$ and $\ b\geq b_{0},\left( ii\right) 0<\protect%
\delta \leq \protect\delta _{0},$ and $\left( iii\right) $every subcube of $%
E $ of side $\protect\delta $ contains a point of $X.$ Here, $\protect\alpha %
_{n}$ deontes the volume of the unit ball in $R^{n}.$}

\subparagraph{The number $c$ is equal to $\frac{b_{0}}{8\protect\gamma _{n}}$
where $\protect\gamma _{n}$ is defined in LEMMA1.The number $C$ is equal to $%
\left( 3^{\frac{3}{4}}\cdot e\cdot \protect\sqrt{2\protect\rho \protect\beta 
}\cdot \protect\sqrt{n}\cdot e^{2n\protect\gamma _{n}}\right) ^{4}\cdot
b_{0}^{\ 3}\cdot \protect\gamma _{n}$, where $\protect\rho =\frac{\protect%
\sqrt{3}}{e}.$ Moreover, $\protect\delta _{0}$ can be defined by}

\begin{center}
$\delta _{0}=\min \left\{ \frac{1}{(3^{\frac{3}{4}}\cdot e\cdot \sqrt{2\rho
\beta }\cdot \sqrt{n}\cdot e^{2n\gamma _{n}})^{4}\cdot b_{0}{}^{3}\cdot
\gamma _{n}},\delta _{n}\right\} $
\end{center}

\subparagraph{where}

\begin{center}
$\delta _{n}=\left\{ 
\begin{array}{c}
\frac{b_{0}}{2\gamma _{n}}\ ~\ \text{if }n=1,\text{ \ \ \ \ \ \ \ \ \ \ \ \
\ \ \ } \\ 
\frac{b_{0}}{2\gamma _{n}(n-1)}\ \text{if }n\text{ is odd, and }n>1, \\ 
\frac{b_{0}}{2\gamma _{n}(n-2)}\ \ \text{if }n\text{ is even, and }n>2, \\ 
\frac{b_{0}}{2\gamma _{n}}\ \ \text{if }n\text{ }=2.\text{ \ \ \ \ \ \ \ \ \
\ \ \ \ }%
\end{array}%
\right\} $
\end{center}

\begin{proof}
\bigskip Let $\delta _{0}$ be as in the theorem. If $\delta \leq \delta
_{0}, $ then$\delta \leq \frac{b_{0}}{2\gamma _{n}}$and$\frac{\gamma
_{n}\delta }{\frac{b_{0}}{2}}\leq 1.$ There is an integer $k\geq 1$ such
that $1\leq (\frac{\gamma _{n}\delta }{\frac{b_{0}}{2}})k\leq 2.\ $This
implies
\end{proof}

\begin{center}
$\frac{b_{0}}{2}\leq \gamma _{n}\delta k\leq b_{0}.$
\end{center}

\subparagraph{Now, let $x$ be any point of the cube $E$ and recall that
Theorem4.2 of $\left[ 8\right] $ implies that}

\begin{flushright}
$\left\vert f(x)-s(x)\right\vert \leq c_{k}\cdot \left\Vert f\right\Vert
_{h}\cdot \tint\limits_{R^{n}}\left\vert y-x\right\vert ^{k}d\left\vert
\sigma \right\vert (y)\ \ \ \ \ \ \ \ \ \ \ \ \ \ \ \ \ \ \ \left( 6\right) $
\end{flushright}

\subparagraph{whenever $k>0,$ where $\protect\sigma $ is any measure
supported on $X$ such that}

\begin{flushright}
$\int p(y)d\sigma (y)=p(x)$ $\ \ \ \ \ \ \ \ \ \ \ \ \ \ \ \ \ \ \ \ \ \ \ \
\ \ \ \ \ \ \ \ \ \ \ \ \ \ \left( 7\right) $
\end{flushright}

\subparagraph{for all polynomials $p$ in $\mathcal{P}_{k-1}.$ Here}

\begin{center}
$c_{k}=\left\{ \tint\limits_{R^{n}}\frac{\left\vert \xi \right\vert ^{2k}}{%
(k!)^{2}}d\mu (\xi )\right\} ^{\frac{1}{2}}$\ \ \ 
\end{center}

\subparagraph{\protect\bigskip}

\subparagraph{whenever $k>0.$}

\subparagraph{LEMMA5. now applies.}

\subparagraph{For odd $n,$}

\subparagraph{$c_{k}=\frac{1}{k!}\left\{ \protect\tint\limits_{R^{n}}\left%
\vert \protect\xi \right\vert ^{2k}d\protect\mu (\protect\xi )\right\} ^{%
\frac{1}{2}}$}

\subparagraph{$\leq \frac{1}{k!}\left\{ \protect\pi ^{\frac{n+1}{2}}\cdot
n\cdot \protect\alpha _{n}\cdot 2^{\frac{2k+n+2}{2}}\cdot \protect\rho ^{%
\frac{2k+n-1}{2}}\cdot \protect\beta ^{k}\cdot (2k+n-1)^{\frac{2k+n-3}{2}%
}\cdot (2+\frac{1}{e})\right\} ^{\frac{1}{2}}$}

\subparagraph{$\leq \frac{1}{k!}\left\{ \protect\pi ^{\frac{n+1}{4}}\cdot
(n\cdot \protect\alpha _{n})^{\frac{1}{2}}\cdot 2^{\frac{n+2}{4}}\cdot 
\protect\rho ^{\frac{n-1}{4}}\cdot (\protect\sqrt{2\protect\rho \protect%
\beta })^{k}\cdot (2k+n-1)^{\frac{2k+n-3}{4}}\cdot \protect\sqrt{2+\frac{1}{e%
}}\right\} $}

\subparagraph{$\leq \frac{\protect\pi ^{\frac{n+1}{4}}\cdot (n\cdot \protect%
\alpha _{n})^{\frac{1}{2}}\cdot 2^{\frac{n+2}{4}}\cdot \protect\rho ^{\frac{%
n-1}{4}}\cdot (\protect\sqrt{2\protect\rho \protect\beta })^{k}\cdot
(2k+n-1)^{\frac{2k+n-3}{4}}\cdot \protect\sqrt{2+\frac{1}{e}}}{\protect\sqrt{%
2\protect\pi }\cdot \protect\rho _{1}^{\ k}\cdot k^{k}}$ by LEMMA 3}

\subparagraph{$=\protect\pi ^{\frac{n-1}{4}}\cdot (n\cdot \protect\alpha %
_{n})^{\frac{1}{2}}\cdot 2^{\frac{n}{4}}\cdot \protect\rho ^{\frac{n-1}{4}%
}\cdot \protect\rho _{3}{}^{k}\cdot (\frac{1}{k^{k}})\cdot (2k+n-1)^{\frac{%
2k+n-3}{4}}\cdot (\protect\sqrt{2+\frac{1}{e}})$}

\subparagraph{where $\protect\rho _{3}=\frac{\protect\sqrt{2\protect\rho 
\protect\beta }}{\protect\rho _{1}}$}

\subparagraph{$=\triangle ^{^{\prime \prime }}\cdot \protect\rho _{3}^{\
k}\cdot \frac{1}{k^{k}}\cdot (2k+n-1)^{\frac{2k+n-3}{4}}$}

\subparagraph{where$\ \triangle ^{^{\prime \prime }}=\protect\sqrt{2+\frac{1%
}{e}}\cdot \protect\pi ^{\frac{n-1}{4}}\cdot (n\cdot \protect\alpha _{n})^{%
\frac{1}{2}}\cdot 2^{\frac{n}{4}}\cdot \protect\rho ^{\frac{n-1}{4}}.$}

\subparagraph{To obtain the desired bound on $\left\vert
f(x)-s(x)\right\vert ,$ it suffices to find a suitable bound for}

\begin{center}
$I=c_{k}\tint\limits_{R^{n}}\left\vert y-x\right\vert ^{k}d\left\vert \sigma
\right\vert (y).$
\end{center}

\subparagraph{Let $Q$ be any cube which contains $x,$ has side $\protect%
\gamma _{n}k\protect\delta ,$ and is contained in $E.$ Subdivide $Q$ into $%
\left( \protect\gamma _{n}k\right) ^{n}$ congruent subcubes of side $\protect%
\delta .$ Since each of these subcubes must contain a point of $X,$ select a
point of $X$ from each subcube and call the resulting discrete set $Y.$ By
Lemma1 we may conclude that there is a measure $\protect\sigma $ supported
on $Y$ which satis fies $\left( 6\right) $ and enjoys the estimate}

\begin{flushright}
$\tint\limits_{R^{n}}d\left\vert \sigma \right\vert (y)\leq e^{2n\gamma
_{n}k}.$ \ \ \ \ \ \ \ \ \ \ \ \ \ \ $\ \ \ \ \ \ \ \ \ \ \ \ \ \ \ \ \ \ \
\ \ \ \left( 8\right) $
\end{flushright}

\subparagraph{We use this measure in $\left( 6\right) $ to obtain an
estimate on $I.$}

\subparagraph{Using $\left( 8\right) $ and the fact that support of $\protect%
\sigma $ is contained in $Q$ whose diameter is $\protect\sqrt{n}\protect%
\gamma _{n}k\protect\delta $ we may write}

\subparagraph{$I\leq \triangle ^{^{\prime \prime }}\cdot \protect\rho %
_{3}^{\ k}\cdot \frac{1}{k^{k}}\cdot (2k+n-1)^{\frac{2k+n-3}{4}}\cdot (%
\protect\sqrt{n}\protect\gamma _{n}k\protect\delta )^{k}\cdot e^{2n\protect%
\gamma _{n}k}$}

\subparagraph{ \ $=\triangle ^{^{\prime \prime }}\cdot (B^{\prime }\protect%
\gamma _{n}\protect\delta )^{k}\cdot (2k+n-1)^{\frac{2k+n-3}{4}}$ where $%
B^{\prime }=\protect\rho _{3}^{\ }\protect\sqrt{n}e^{^{2n\protect\gamma %
_{n}}}$}

\subparagraph{$\ \leq \triangle ^{^{\prime \prime }}\cdot (B^{\prime }%
\protect\gamma _{n}\protect\delta )^{k}\cdot (3k)^{\frac{3k}{4}}$ \ $\left(
if\text{ \ }k\geq n-1\right) $}

\subparagraph{$\ =\triangle ^{^{\prime \prime }}\cdot (B^{\prime }\protect%
\gamma _{n}\protect\delta 3^{\frac{3}{4}}k^{\frac{3}{4}})^{k}$}

\subparagraph{\ \ $\leq \triangle ^{^{\prime \prime }}\cdot \left[ 3^{\frac{3%
}{4}}B^{\prime }\protect\gamma _{n}\protect\delta (\frac{b_{0}}{\protect%
\gamma _{n}\protect\delta })^{\frac{3}{4}}\right] ^{k}$ \ $($since $k\leq 
\frac{b_{0}}{\protect\gamma _{n}\protect\delta })$}

\subparagraph{\ \ $=\triangle ^{^{\prime \prime }}\cdot \left[ B^{^{^{\prime
\prime }}}\protect\gamma _{n}^{~\frac{1}{4}}b_{0}^{\frac{3}{4}}\protect%
\delta ^{\frac{1}{4}}\right] ^{k}$ where $B^{^{\prime \prime }}=3^{\frac{3}{4%
}}B^{\prime }$}

\subparagraph{ \ $=\triangle ^{^{\prime \prime }}\cdot \left[ B^{^{^{\prime
\prime }}}b_{0}^{\frac{3}{4}}\protect\gamma _{n}^{~\frac{1}{4}}\protect%
\delta ^{\frac{1}{4}}\right] ^{\frac{b_{0}}{2\protect\gamma _{n}\protect%
\delta }}$ $($since $\frac{b_{0}}{2\protect\gamma _{n}\protect\delta }\leq
k) $}

\subparagraph{where $B^{^{^{\prime \prime }}}b_{0}^{\frac{3}{4}}\protect%
\gamma _{n}^{~\frac{1}{4}}\protect\delta ^{\frac{1}{4}}\leq 1$ if and only
if $B^{^{\prime \prime ^{4}}}b_{0}^{3}\protect\gamma _{n}\protect\delta \leq
1$ and is guaranteed by $\protect\delta \leq \protect\delta _{0}$}

\subparagraph{\ \ $=\triangle ^{^{\prime \prime }}\cdot \left[ B^{^{^{\prime
\prime ^{4}}}}b_{0}^{3}\protect\gamma _{n}^{~}\protect\delta \right] ^{\frac{%
b_{0}}{8\protect\gamma _{n}\protect\delta }}$}

\subparagraph{ \ $=\triangle ^{^{\prime \prime }}\cdot \left[ C\protect%
\delta \right] ^{\frac{c}{\protect\delta }}$ where $C=B^{^{^{\prime \prime
^{4}}}}b_{0}^{3}\protect\gamma _{n}^{~}$ and $c=\frac{b_{0}}{8\protect\gamma %
_{n}}.$}

\subparagraph{Note that $k\geq n-1$ is guaranteed by $n-1\leq \frac{b_{0}}{2%
\protect\gamma _{n}\protect\delta }\leq k.$}

\subparagraph{For even $n,$}

\subparagraph{$c_{k}=\frac{1}{k!}\left\{ \protect\tint\limits_{R^{n}}\left%
\vert \protect\xi \right\vert ^{2k}d\protect\mu (\protect\xi )\right\} ^{%
\frac{1}{2}}$}

\subparagraph{$\ \ \ \leq \frac{1}{k!}\left\{ \protect\pi ^{\frac{n+1}{2}%
}\cdot n\cdot \protect\alpha _{n}\cdot 2^{\frac{2k+n+3}{2}}\cdot \protect%
\rho ^{\frac{2k+n-2}{2}}\cdot \protect\beta ^{k}\cdot (2k+n-2)^{\frac{2k+n-4%
}{2}}\right\} ^{\frac{1}{2}}$}

\subparagraph{$\ \ \ =\frac{1}{k!}\left\{ \protect\pi ^{\frac{n+1}{4}}\cdot
(n\protect\alpha _{n})^{\frac{1}{2}}\cdot 2^{\frac{n+3}{4}}\cdot \protect%
\rho ^{\frac{n-2}{4}}\cdot (\protect\sqrt{2\protect\rho \protect\beta }%
)^{k}\cdot (2k+n-2)^{\frac{2k+n-4}{4}}\right\} $}

\subparagraph{$\protect\vspace{1pt}$}

\subparagraph{$\ \ \ \leq \frac{\protect\pi ^{\frac{n+1}{4}}\cdot (n\protect%
\alpha _{n})^{\frac{1}{2}}\cdot 2^{\frac{n+3}{4}}\cdot \protect\rho ^{\frac{%
n-2}{4}}\cdot (\protect\sqrt{2\protect\rho \protect\beta })^{k}\cdot
(2k+n-2)^{\frac{2k+n-4}{4}}}{\protect\sqrt{2\protect\pi }\cdot \protect\rho %
_{1}^{\ k}\cdot k^{k}}$}

\subparagraph{$\ \ \ =\triangle ^{^{\prime \prime }}\cdot \protect\rho %
_{3}^{\ k}\cdot \frac{(2k+n-2)^{\frac{2k+n-4}{4}}}{k^{k}}$ where $\triangle
^{^{\prime \prime }}=\protect\pi ^{\frac{n-1}{4}}\cdot (n\protect\alpha %
_{n})^{\frac{1}{2}}\cdot 2^{\frac{n+1}{4}}\cdot \protect\rho ^{\frac{n-2}{4}%
} $}

\subparagraph{$\ $and $\protect\rho _{3}=\frac{\protect\sqrt{2\protect\rho 
\protect\beta }}{\protect\rho _{1}}$}

\subparagraph{Now,}

\subparagraph{$\ \ \ \ \ \ \ \ \ \ \ \ \ \ \ \ \ \ \ \ \ \ I=c_{k}\protect%
\tint\limits_{R^{n}}\left\vert y-x\right\vert ^{k}d\left\vert \protect\sigma %
\right\vert (y)$}

\subparagraph{$\ \ \ \ \ \ \ \ \ \ \ \ \ \ \ \ \ \ \ \ \ \ \ \leq \triangle
^{^{\prime \prime }}\cdot \protect\rho _{3}^{\ k}\cdot \frac{(2k+n-2)^{\frac{%
2k+n-4}{4}}}{k^{k}}\cdot (\protect\sqrt{n}\cdot \protect\gamma _{n}\cdot
k\cdot \protect\delta )^{k}\cdot e^{2n\protect\gamma _{n}k}$}

\subparagraph{$\ \ \ \ \ \ \ \ \ \ \ \ \ \ \ \ \ \ \ \ \ \ \ =\triangle
^{^{\prime \prime }}\cdot (B^{\prime }\protect\gamma _{n}\protect\delta %
)^{k}\cdot (2k+n-2)^{\frac{2k+n-4}{4}}$ where $B^{\prime }=\protect\rho %
_{3}^{\ }\protect\sqrt{n}e^{^{2n\protect\gamma _{n}}}$}

\subparagraph{$\ \ \ \ \ \ \ \ \ \ \ \ \ \ \ \ \ \ \ \ \ \ \ \leq \triangle
^{^{\prime \prime }}\cdot (B^{\prime }\protect\gamma _{n}\protect\delta %
)^{k}\cdot (3k)^{\frac{3k}{4}}$ \ $\left( whenever\text{ \ }k\geq n-2\right) 
$}

\subparagraph{$\ \ \ \ \ \ \ \ \ \ \ \ \ \ \ \ \ \ \ \ \ \ \ =\triangle
^{^{\prime \prime }}\cdot (B^{\prime }\protect\gamma _{n}\protect\delta 3^{%
\frac{3}{4}}k^{\frac{3}{4}})^{k}$}

\subparagraph{\ \ \ \ \ \ \ \ \ \ \ \ \ \ \ \ \ \ \ \ $\ \leq \triangle
^{^{\prime \prime }}\cdot \left[ 3^{\frac{3}{4}}B^{\prime }\protect\gamma %
_{n}\protect\delta (\frac{b_{0}}{\protect\gamma _{n}\protect\delta })^{\frac{%
3}{4}}\right] ^{k}$ \ $(k\leq \frac{b_{0}}{\protect\gamma _{n}\protect\delta 
})$}

\subparagraph{\ \ \ \ \ \ \ \ \ \ \ \ \ \ \ \ \ \ \ \ \ \ $=\triangle
^{^{\prime \prime }}\cdot \left[ B^{^{^{\prime \prime }}}\protect\gamma %
_{n}^{~\frac{1}{4}}b_{0}^{\frac{3}{4}}\protect\delta ^{\frac{1}{4}}\right]
^{k}$ where $B^{^{\prime \prime }}=3^{\frac{3}{4}}\cdot B^{\prime }$}

\subparagraph{\ \ \ \ \ \ \ \ \ \ \ \ \ \ \ \ \ \ \ \ \ \ $\leq ~\triangle
^{^{\prime \prime }}\cdot \left[ B^{^{\prime \prime }}b_{0}^{\frac{3}{4}}%
\protect\gamma _{n}^{~\frac{1}{4}}\protect\delta ^{\frac{1}{4}}\right] ^{%
\frac{b_{0}}{2\protect\gamma _{n}\protect\delta }}$ since $\frac{b_{0}}{2%
\protect\gamma _{n}\protect\delta }\leq k$}

\subparagraph{ \ \ \ \ \ \ \ \ \ \ \ \ \ \ \ \ \ \ \ \ \ \ \ \ \ and $%
B^{^{\prime \prime ^{4}}}b_{0}^{\frac{3}{4}}\protect\gamma _{n}^{~\frac{1}{4}%
}\protect\delta ^{\frac{1}{4}}\leq 1$ due to $\protect\delta \leq \protect%
\delta _{0}.$}

\subparagraph{This gives}

\subparagraph{$\ \ \ \ \ \ \ \ \ \ \ \ \ \ \ \ \ \ \ \ \ \ \ \ I\leq
\triangle ^{^{\prime \prime }}\left[ B^{^{^{\prime \prime ^{4}}}}b_{0}^{~3}%
\protect\gamma _{n}^{~}\protect\delta \right] ^{\frac{b_{0}}{8\protect\gamma %
_{n}\protect\delta }}$}

\subparagraph{$\ \ \ \ \ \ \ \ \ \ \ \ \ \ \ \ \ \ \ \ \ \ \ \ \ =\triangle
^{^{\prime \prime }}\cdot \left[ C\protect\delta \right] ^{\frac{c}{\protect%
\delta }}$ where $C=B^{^{^{\prime \prime ^{4}}}}b_{0}^{3}\protect\gamma %
_{n}^{~}$ and $c=\frac{b_{0}}{8\protect\gamma _{n}}.$}

\subparagraph{Note that $k\geq n-2$ is guaranteed by $\protect\delta \leq 
\protect\delta _{0}\leq \frac{b_{0}}{2\protect\gamma _{n}(n-2)}$and $n-2\leq 
\frac{b_{0}}{2\protect\gamma _{n}\protect\delta }\leq k.$}

\subparagraph{We conclude that}

\begin{center}
$\left\vert f(x)-s(x)\right\vert $ $\leq \triangle ^{^{\prime \prime }}\cdot %
\left[ C\delta \right] ^{\frac{c}{\delta }}\left\Vert f\right\Vert _{h}$
\end{center}

\subparagraph{whenever $0<$ $\protect\delta \leq \protect\delta _{0}$ as
stated in the theorem.}

\subparagraph{Remark}

\textit{The high-level error bound for Gaussians was first put forward by
Madych and Nelson in Theorem 3. of }$\left[ 9\right] .$\textit{\ However
their proof is incomplete. In their theorem the radial basis function }$h$%
\textit{\ must satisfy the key condition}

\begin{flushright}
$\tint\limits_{R^{n}}\left\vert \xi \right\vert ^{k}d\mu (\xi )\leq \rho
^{k}k^{rk}$ $\ \ \ \ \ \ \ \ \ \ \ \ \ \ \ \ \ \ \ \ \ \ \ \ \ \ \ \ \ \ \ \
\ \ \ \ (9)$
\end{flushright}

\subparagraph{for $k>2m,$ where $r$ is a real constant and $\protect\rho $
is a positive constant. As pointed out by them in page 102 of $\left[ 9%
\right] ,$ $r=\frac{1}{2}$ if $h$ is a Gaussian function. Moreover, they
only treated the case $\protect\beta =1.$ However$\left( 9\right) $ holds
only when $n=1.$ For $n>1,\left( 9\right) $ should be replaced by Lemma 5 of
this paper. Consequently their high-level error bound is essentially only
suitable for $R^{1}.$}

\subparagraph{Note that Theorem 1 can be improved if $n$ is known in
advance, especially when $n=1$ or $2.$ For example, suppose $n=1.$ Let}

\begin{center}
\bigskip $\delta _{0}=\min \left\{ \frac{1}{(\sqrt{2}\cdot e\cdot \sqrt{%
2\rho \beta }\cdot \sqrt{n}\cdot e^{2n\gamma _{n}})^{2}\cdot {}\gamma
_{n}\cdot b_{0}},\frac{b_{0}}{2\gamma _{n}}\right\} .$
\end{center}

\subparagraph{Then for $\protect\delta \leq \protect\delta _{0},$}

\begin{center}
$\ \ \ \ \ \ \ \ \ \ \ \ \ \ \ \ \ \ \ \ I\leq \triangle ^{^{\prime \prime
}}\cdot (B^{\prime }\gamma _{n}\delta )^{k}\cdot (2k)^{\frac{2k-2}{4}}$ \ 

$\ \ \ \ \ \ \ \ \ \ \ \ \ \ \ \ =\triangle ^{^{\prime \prime }}\cdot
(B^{\prime }\gamma _{n}\delta )^{k}\cdot (2k)^{\frac{k}{2}}$

\ \ \ \ \ \ \ \ \ \ \ \ \ \ \ \ \ \ \ \ $\leq \triangle ^{^{\prime \prime
}}\cdot \left[ B^{\prime }\gamma _{n}\delta \cdot \sqrt{2}\sqrt{k}\right]
^{k}$

\ \ \ \ \ \ \ \ \ \ \ \ \ \ \ \ \ \ \ \ \ \ \ \ $=\triangle ^{^{\prime
\prime }}\cdot \left[ \sqrt{2}B^{\prime }\gamma _{n}^{~}\delta (\frac{b_{0}}{%
\gamma _{n}\delta })^{\frac{1}{2}}\right] ^{k}$

\ \ \ \ \ \ \ \ \ \ \ \ \ \ \ \ \ \ \ \ \ \ \ \ \ \ \ \ \ \ \ \ \ \ \ \ \ \
\ \ \ \ \ \ \ $=\triangle ^{^{^{\prime \prime }}}\cdot \left[ B^{^{\prime
\prime }}\gamma _{n}^{~\frac{1}{2}}\delta ^{\frac{1}{2}}b_{0}^{\frac{1}{2}}%
\right] ^{k}$ where $B^{^{\prime \prime }}=\sqrt{2}B^{^{\prime }}$

\ \ \ \ \ \ \ \ \ \ \ \ \ \ \ \ \ \ \ \ \ \ \ \ $=\triangle ^{^{^{\prime
\prime }}}\cdot \left[ B^{^{\prime \prime }}\gamma _{n}^{~\frac{1}{2}}\delta
^{\frac{1}{2}}b_{0}^{\frac{1}{2}}\right] ^{\frac{b_{0}}{2\gamma _{n}\delta }%
} $ \ \ 
\end{center}

\subparagraph{where $B^{^{\prime \prime }}\protect\gamma _{n}^{~\frac{1}{2}}%
\protect\delta ^{\frac{1}{2}}b_{0}^{\frac{1}{2}}\leq 1$ is guaranteed by $%
\protect\delta \leq \protect\delta _{0},$}

\begin{center}
\ \ \ \ \ \ \ \ \ \ \ \ \ \ \ \ \ \ \ \ \ \ \ \ \ \ \ \ \ \ \ \ \ \ \ \ \ \
\ \ $=~\triangle ^{^{^{\prime \prime }}}\cdot \left[ B^{^{^{\prime \prime
}}}\gamma _{n}^{~}b_{0}\delta \right] ^{\frac{b_{0}}{4\gamma _{n}\delta }%
}=\triangle ^{^{^{\prime \prime }}}\cdot \left[ C\delta \right] ^{\frac{c}{%
\delta }}$
\end{center}

\subparagraph{where $C=B^{^{^{\prime \prime ^{2}}}}\protect\gamma %
_{n}^{~}b_{0}$ and $c=\frac{b_{0}}{4\protect\gamma _{n}}.$}

\subparagraph{In this case, $\protect\delta _{0},C,c$ are the same as Madych
and Nelson's results without any sacrifice.}

\subparagraph{What's noteworthy is that in Theorem 1 the parameter $\protect%
\delta $ is not the generally used fill distance. For easy use we should
transform Theorem 1 into a statement described by the fill distance.}

\subparagraph{Let}

\begin{center}
\ \ \ \ \ $d(\Omega ,X)=\underset{y\in \Omega }{\sup }$ $\underset{x\in X}{%
\inf }\left\vert y-x\right\vert $
\end{center}

\subparagraph{be the fill distane. Observe that every cube of side $\protect%
\delta $ contains a ball of radius $\frac{\protect\delta }{2}.$ Thus the
subcube condition in Theorem1 is satisfied when $\protect\delta =2d(E,X).$
More generally, we can easily conclude the following:}

\subparagraph{\textbf{Corollary}}

\textit{Suppose that }$h$\textit{\ satisfies the hypotheses of the theorem, }%
$\Omega $\textit{\ is a set which can be expressed as the union of rotations
and translations of a fixed cube of side }$b_{0},$\textit{\ and }$X$\textit{%
\ is a subset of }$R^{n}.$\textit{\ Then there are positive contants }$%
d_{0},~c^{\prime }~,$\textit{and }$C^{\prime }$\textit{\ for which the
following is true:}\bigskip {}

\textbf{If }$f\in C_{h,m}$\textbf{\ and }$s$\textbf{\ is the }$h$\textbf{\
spline that interpolates }$f$\textbf{\ on a subset }$X$\textbf{\ of }$R^{n},$%
\textbf{then}

\begin{center}
$\left\vert f(x)-s(x)\right\vert $ $\leq \triangle ^{^{\prime \prime }}\left[
C^{\prime }d\right] ^{\frac{c^{^{\prime }}}{d}}\left\Vert f\right\Vert _{h}$
\end{center}

\subparagraph{,where $\triangle ^{^{\prime \prime }}$ is as in Theorem1 and $%
C^{\prime }=2C,c^{\prime }=\frac{c}{2}$ with $C$ and $c$ defined in
Theorem1, holds for all $x$ in a cube $E\subseteq \Omega $ provided that$%
\left( i\right) $ $E$ has side $b$ and $b\geq b_{0},\left( ii\right) 0<d\leq
d_{0},$and $\left( iii\right) $every subcube of $E$ of side $2d$ contains a
point of $X.$}

\textbf{Here }$d$\textbf{\ denotes }$d(\Omega ,X)$\textbf{\ or }$d(E,X)$%
\textbf{\ and }$m=0$\textbf{\ since Gaussians are c.p.d of order }$0.$

\begin{proof}
Let $d_{0}=\frac{\delta _{0}}{2}$ and $\delta =2d$ where $\delta _{0}$ is as
in the theorem. If $0<d\leq d_{0},0<$ $\delta \leq 2d_{0}=\delta _{0},$ By
the theorem,

$\ \ \ \ \ \ \ \ \ \ \ \ \ \ \ \ \ \ \ \left\vert f(x)-s(x)\right\vert $ $%
\leq \triangle ^{^{\prime \prime }}\left[ C\delta \right] ^{\frac{c}{\delta }%
}\cdot \left\Vert f\right\Vert _{h}$

\ \ \ \ \ \ \ \ \ \ \ \ \ \ \ \ \ $\ \ \ \ \ \ \ \ \ \ \ \ \ \ \ \ \ \
=\triangle ^{^{\prime \prime }}\left[ C2d\right] ^{\frac{c}{2d}}\cdot
\left\Vert f\right\Vert _{h}$

\ \ \ \ \ \ \ \ \ \ \ \ \ \ \ \ \ \ \ \ \ \ \ \ \ \ \ \ \ \ \ \ \ \ $\
=\triangle ^{^{\prime \prime }}\left[ C^{\prime }d\right] ^{\frac{%
c^{^{\prime }}}{\delta }}\cdot \left\Vert f\right\Vert _{h}$

holds for all $x$ in $\Omega .\ $ $\ \ \ \ \ \ \ \ \ \ \ \ \ \ \ \ \ \ \ \ \
\ \ \ \ \ \ \ \ \ \ \ \ \ \ \ \ \ \ \ \ \ \ \ \ \ \ \ \ \ \ \ \ \ \ \ \ \ \
\ \ \ \ \ \ \ \ \ \ \ \ \ \ \ \ \ \ \ \ \ \ \ \ \ \ \ \ \ \ \ \ \ \ \ \ \ \
\ \ \ \ \ \ \ \ \ \ \ \ \ \ \ \ \ \ \ \ \ \ \ \ \ \ \ \ \ \ \ \ \ \ \ \ \ \
\ \ \ \ \ \ \ \ \ \ $
\end{proof}

\begin{flushright}
$\square $
\end{flushright}

\subparagraph{Remark}

\textit{The space }$C_{h,m}$\textit{\ probably is unfamiliar to most people.
It's introduced by Madych and Nelson in }$\left[ 7\right] $\textit{\ and }$%
\left[ 8\right] .$\textit{Later Luh made characterizations for it in }$\left[
3\right] $\textit{\ and }$\left[ 4\right] .$\textit{\ Many people think that
it's defined by Gelfand and Shilov's definition of generalized Fourier
transform, and is therefore difficult to deal with. This is wrong. In fact,
it can be characterized by Schwartz's definition of generalized Fourier
transform. The situation is not so bad. Moreover, many people mistake }$%
C_{h,m}$\textit{\ to be the closure of Wu and Schaback's function space
which is defined in }$\left[ 10\right] .$\textit{\ This is also wrong. The
two spaces have very little connection. Luh also made a clarification for
this problem. For further details, please see }$\left[ 5\right] $\textit{\
and }$\left[ 6\right] .$

\begin{center}
{\Large REFERENCES}
\end{center}

\begin{enumerate}
\item J. von zur Gathen and J. Gerhard,\textquotedblleft\ Modern Computer
Algebra\textquotedblright

Second Edition, 2003, CAMBRIDGE University Press.

\item R. L. GRAHAM, D.E. KNUTH, and O. PATASHNIK(1994),

\textquotedblleft Concrete Mathematics,\textquotedblright\ Addison-Wesley,
Reading MA, 2nd edition, First edition 1989.{}

\item Lin-Tain Luh, The Equivalence Theory of Native Spaces, Approximations
Theory and its Applications, 2001, 17:1, 76-96.

\item Lin-Tain Luh, The Embedding Theory of Native Spaces, Approximations
Theory and its Applications, 2001, 17:4, 90-104.

\item Lin-Tain Luh, On Wu and Schaback's Error Bound, to appear.

\item Lin-Tain Luh, The Completeness of Function Spaces, to appear.

\item W.R.MADYCH AND S.A.NELSON, Multivariate interpolation and
conditionally positive definite function, Approx. Theory Appl. 4, No.
4(1988), 77-89.

\item W.R.MADYCH AND S.A.NELSON, Multivariate interpolation and
conditionally positive definite function, II, Math. Comp. 54(1990), 211-230.

\item W.R.MADYCH AND S.A.NELSON, Bounds on Multivariate Polynomials and
Exponential Error Estimates for Multiquadric Interpolation, J. Approx.
Theory 70, 1992, 94-114.

\item Z. Wu and R. Schaback, Local Error Estimates for Radial Basis Function
Interpolation of Scattered Data, IMA J. of Numerical Analysis, 1993,
13:13-27.
\end{enumerate}

\end{document}